\documentclass{amsart}
\usepackage{amsmath,amssymb}
\usepackage{mathpazo}
\usepackage{ifthen}
\usepackage{tikz,pgf}
    \usetikzlibrary{shapes,calc,decorations.pathmorphing}
	\pgfdeclarelayer{background}
	\pgfdeclarelayer{middle}
	\pgfsetlayers{background,middle,main}
\usepackage[pagebackref,
    colorlinks=true,
	urlcolor=blue!20!black,
	citecolor=green!50!black]{hyperref}
\usepackage[nobysame]{amsrefs}
\newtheorem{theorem}{Theorem}[section]
\newtheorem{proposition}[theorem]{Proposition}

\newtheorem{lemma}[theorem]{Lemma}

\theoremstyle{remark}

\newtheorem{example}[theorem]{Example}

\newtheorem{construction}[theorem]{Construction}
\newcommand{\defn}[1]{{\color{green!50!black}\emph{#1}}}
\newcommand{\defs}{\stackrel{\mathsf{def}}{=}}
\newcommand{\wb}{\mathfrak{w}}
\newcommand{\Words}{\mathsf{W}}
\newcommand{\Tilings}{\mathsf{T}}

\newcommand{\Tile}[3]{
    \filldraw[draw=white!50!black,fill=#3]({#1*\s},0*\s) -- ({(#1+#2)*\s},0*\s) -- ({(#1+#2)*\s},1*\s) -- ({#1*\s},1*\s) -- cycle;
    \foreach \a in {1,...,#2}{
        \ifthenelse{\equal{\a}{#2}}{}{
            \draw[dotted,white!50!black]({(\a+#1)*\s},0*\s) -- ({(\a+#1)*\s},1*\s);
        }
    }
}

\makeatletter
\let\orgdescriptionlabel\descriptionlabel
\renewcommand*{\descriptionlabel}[1]{%
  \let\orglabel\label
  \let\label\@gobble
  \phantomsection
  \protected@edef\@currentlabel{#1\unskip}%
  \let\label\orglabel
  \orgdescriptionlabel{(#1)}%
}
\makeatother

\title{On Two-Toned Tilings and $(m,n)$-Words}
\author{Henri M{\"u}hle}
\address{KLA Dresden City Centre, Dr.-K{\"u}lz-Ring 15, 01067 Dresden, Germany}
\email{henri.muehle@proton.me}

\begin{document}

\begin{abstract}
    In this article, we describe an explicit bijection between the set of $(m,n)$-words as defined by Pilaud and Poliakova and the set of of two-toned tilings of a strip of length $m+n$. 
\end{abstract}

\maketitle

\section{Introduction}
For two integers $m$ and $n$, V.~Pilaud and D.~Poliakova introduced so-called $(m,n)$-words as intermediate objects in their definition of Hochschild polytopes \cite{pilaud25hochschild}. These words were counted in \cite{muehle24combinatorics,pilaud25hochschild}, and it was observed computationally by T.~Copeland~\cite{copeland25communication} that the number of $(m,n)$-words agrees with the coefficient of $x^n$ in $\left(\frac{1-x}{1-2x}\right)^{m+1}$.

In \cite{benjamin11combinatorics}, two-toned tilings of a strip were introduced and studied, and it was shown that a particular class of two-toned tilings (namely those using $m$ squares of one color and arbitrary strips of cumulated length $n$ of another color) is enumerated by the coefficients of $\left(\frac{1-x}{1-2x}\right)^{m+1}$.

The main purpose of this article is the explicit construction of a bijection between the set of $(m,n)$-words and the set of two-toned tilings of a strip of length $m+n$.

\section{Basics}
Throughout this article, we use the abbreviation $[k]\defs\{1,2,\ldots,k\}$ for a positive integer $k$.

\subsection{$(m,n)$-Words}
Let $m,n\geq 0$. Following \cite[Definition~77]{pilaud25hochschild}, an \defn{$(m,n)$-word} is a word $w_1w_2\ldots w_n$ of length $n$ over the alphabet $\{0,1,\ldots,m+1\}$ such that 
\begin{description}
    \item[MN1\label{it:mn_1}] $w_1\neq m+1$
    \item[MN2\label{it:mn_2}] for $1\leq s\leq m$, $w_i=s$ implies $w_j\geq s$ for all $j<i$.
\end{description}

In other words, an $(m,n)$-word is a weakly decreasing sequence of length $n$ of numbers in $\{0,1,\ldots,m\}$, where some of the entries, except for the first one, can be replaced by $m+1$. Then, a \defn{topless} $(m,n)$-word is an $(m,n)$-word that does not contain the letter $m+1$. 

\begin{lemma}[{\cite[Proposition~16]{muehle24combinatorics}}]\label{lem:mn_words_cardinality}
    For $m\geq 0$, $n\geq 1$, the number of $(m,n)$-words is
    \begin{displaymath}
        \sum_{k=1}^{n}\binom{m+k}{k}\binom{n-1}{k-1}.
    \end{displaymath}
\end{lemma}

\begin{example}
    The 25 $(2,3)$-words are the following:
    \begin{displaymath}\begin{aligned}
        & 000, && 003, && 030, && 033, && 100, && 103, && 110, && 111, && 113, && 130, && 131, && 133, && 200,\\
        & 203, && 210, && 211, && 213, && 220, && 221, && 222, && 223, && 230, && 231, && 232, && 233.
    \end{aligned}\end{displaymath}
\end{example}

\subsection{Two-toned tilings}

A \defn{strip} of length $k$ is a $1\times k$-rectangle. A strip of length $1$ is a \defn{square}. If $S$ is any strip, then we sometimes use $\lvert S\rvert$ for its length.

A \defn{tiling} of a strip of length $k$ is a collection of strips of lengths $k_1,k_2,\ldots,k_s$ such that $k_1+k_2+\cdots+k_s=k$. 

A \defn{two-toned tiling} of length $m+n$ is a tiling of a strip of length $m+n$ into $m$ red squares and arbitrarily many blue strips. This is to imply that the sum of lengths of the blue strips is $n$. 

Let us denote the set of two-toned tilings of length $m+n$ by $\Tilings(m,n)$. 

\begin{lemma}[{\cite[Identity~3]{benjamin11combinatorics}}]\label{lem:two_toned_tilings_cardinality}
    For $m\geq 0$, $n\geq 1$, the number of two-toned tilings of length $m+n$ is
    \begin{displaymath}
        \sum_{k=1}^{n}\binom{m+k}{k}\binom{n-1}{k-1}.
    \end{displaymath}
\end{lemma}

\begin{lemma}[{\cite[Equation~2.1]{davis20further}}]
    For $m\geq 0$ and $n\geq 1$, the number of two-toned tilings of $m+n$ is the coefficient of $x^{n}$ in $\left(\frac{1-x}{1-2x}\right)^{m+1}$.
\end{lemma}

\begin{example}
    The $25$ two-toned tilings of length $2+3$ are the following:
    \begin{displaymath}\begin{aligned}
        & 
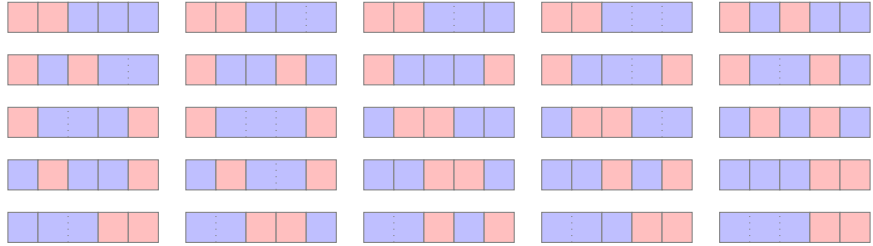
\begin{tikzpicture} 
            \def\s{0.4};
            \Tile{0}{1}{red!25!white}
            \Tile{1}{1}{red!25!white}
            \Tile{2}{1}{blue!25!white}
            \Tile{3}{1}{blue!25!white}
            \Tile{4}{1}{blue!25!white}
        \end{tikzpicture}
        && \begin{tikzpicture} 
            \def\s{0.4};
            \Tile{0}{1}{red!25!white}
            \Tile{1}{1}{red!25!white}
            \Tile{2}{1}{blue!25!white}
            \Tile{3}{2}{blue!25!white}
        \end{tikzpicture}
        && \begin{tikzpicture} 
            \def\s{0.4};
            \Tile{0}{1}{red!25!white}
            \Tile{1}{1}{red!25!white}
            \Tile{2}{2}{blue!25!white}
            \Tile{4}{1}{blue!25!white}
        \end{tikzpicture}
        && \begin{tikzpicture} 
            \def\s{0.4};
            \Tile{0}{1}{red!25!white}
            \Tile{1}{1}{red!25!white}
            \Tile{2}{3}{blue!25!white}
        \end{tikzpicture}
        && \begin{tikzpicture} 
            \def\s{0.4};
            \Tile{0}{1}{red!25!white}
            \Tile{1}{1}{blue!25!white}
            \Tile{2}{1}{red!25!white}
            \Tile{3}{1}{blue!25!white}
            \Tile{4}{1}{blue!25!white}
        \end{tikzpicture}\\
        & \begin{tikzpicture} 
            \def\s{0.4};
            \Tile{0}{1}{red!25!white}
            \Tile{1}{1}{blue!25!white}
            \Tile{2}{1}{red!25!white}
            \Tile{3}{2}{blue!25!white}
        \end{tikzpicture}
        && \begin{tikzpicture} 
            \def\s{0.4};
            \Tile{0}{1}{red!25!white}
            \Tile{1}{1}{blue!25!white}
            \Tile{2}{1}{blue!25!white}
            \Tile{3}{1}{red!25!white}
            \Tile{4}{1}{blue!25!white}
        \end{tikzpicture}
        && \begin{tikzpicture} 
            \def\s{0.4};
            \Tile{0}{1}{red!25!white}
            \Tile{1}{1}{blue!25!white}
            \Tile{2}{1}{blue!25!white}
            \Tile{3}{1}{blue!25!white}
            \Tile{4}{1}{red!25!white}
        \end{tikzpicture}
        && \begin{tikzpicture} 
            \def\s{0.4};
            \Tile{0}{1}{red!25!white}
            \Tile{1}{1}{blue!25!white}
            \Tile{2}{2}{blue!25!white}
            \Tile{4}{1}{red!25!white}
        \end{tikzpicture}
        && \begin{tikzpicture} 
            \def\s{0.4};
            \Tile{0}{1}{red!25!white}
            \Tile{1}{2}{blue!25!white}
            \Tile{3}{1}{red!25!white}
            \Tile{4}{1}{blue!25!white}
        \end{tikzpicture}\\
        & \begin{tikzpicture} 
            \def\s{0.4};
            \Tile{0}{1}{red!25!white}
            \Tile{1}{2}{blue!25!white}
            \Tile{3}{1}{blue!25!white}
            \Tile{4}{1}{red!25!white}
        \end{tikzpicture}
        && \begin{tikzpicture} 
            \def\s{0.4};
            \Tile{0}{1}{red!25!white}
            \Tile{1}{3}{blue!25!white}
            \Tile{4}{1}{red!25!white}
        \end{tikzpicture}
        && \begin{tikzpicture} 
            \def\s{0.4};
            \Tile{0}{1}{blue!25!white}
            \Tile{1}{1}{red!25!white}
            \Tile{2}{1}{red!25!white}
            \Tile{3}{1}{blue!25!white}
            \Tile{4}{1}{blue!25!white}
        \end{tikzpicture}
        && \begin{tikzpicture} 
            \def\s{0.4};
            \Tile{0}{1}{blue!25!white}
            \Tile{1}{1}{red!25!white}
            \Tile{2}{1}{red!25!white}
            \Tile{3}{2}{blue!25!white}
        \end{tikzpicture}
        && \begin{tikzpicture} 
            \def\s{0.4};
            \Tile{0}{1}{blue!25!white}
            \Tile{1}{1}{red!25!white}
            \Tile{2}{1}{blue!25!white}
            \Tile{3}{1}{red!25!white}
            \Tile{4}{1}{blue!25!white}
        \end{tikzpicture}\\
        & \begin{tikzpicture} 
            \def\s{0.4};
            \Tile{0}{1}{blue!25!white}
            \Tile{1}{1}{red!25!white}
            \Tile{2}{1}{blue!25!white}
            \Tile{3}{1}{blue!25!white}
            \Tile{4}{1}{red!25!white}
        \end{tikzpicture} 
        && \begin{tikzpicture} 
            \def\s{0.4};
            \Tile{0}{1}{blue!25!white}
            \Tile{1}{1}{red!25!white}
            \Tile{2}{2}{blue!25!white}
            \Tile{4}{1}{red!25!white}
        \end{tikzpicture}
        && \begin{tikzpicture} 
            \def\s{0.4};
            \Tile{0}{1}{blue!25!white}
            \Tile{1}{1}{blue!25!white}
            \Tile{2}{1}{red!25!white}
            \Tile{3}{1}{red!25!white}
            \Tile{4}{1}{blue!25!white}
        \end{tikzpicture}
        && \begin{tikzpicture} 
            \def\s{0.4};
            \Tile{0}{1}{blue!25!white}
            \Tile{1}{1}{blue!25!white}
            \Tile{2}{1}{red!25!white}
            \Tile{3}{1}{blue!25!white}
            \Tile{4}{1}{red!25!white}
        \end{tikzpicture}
        && \begin{tikzpicture} 
            \def\s{0.4};
            \Tile{0}{1}{blue!25!white}
            \Tile{1}{1}{blue!25!white}
            \Tile{2}{1}{blue!25!white}
            \Tile{3}{1}{red!25!white}
            \Tile{4}{1}{red!25!white}
        \end{tikzpicture}\\
        & \begin{tikzpicture} 
            \def\s{0.4};
            \Tile{0}{1}{blue!25!white}
            \Tile{1}{2}{blue!25!white}
            \Tile{3}{1}{red!25!white}
            \Tile{4}{1}{red!25!white}
        \end{tikzpicture}
        && \begin{tikzpicture} 
            \def\s{0.4};
            \Tile{0}{2}{blue!25!white}
            \Tile{2}{1}{red!25!white}
            \Tile{3}{1}{red!25!white}
            \Tile{4}{1}{blue!25!white}
        \end{tikzpicture}
        && \begin{tikzpicture} 
            \def\s{0.4};
            \Tile{0}{2}{blue!25!white}
            \Tile{2}{1}{red!25!white}
            \Tile{3}{1}{blue!25!white}
            \Tile{4}{1}{red!25!white}
        \end{tikzpicture}
        && \begin{tikzpicture} 
            \def\s{0.4};
            \Tile{0}{2}{blue!25!white}
            \Tile{2}{1}{blue!25!white}
            \Tile{3}{1}{red!25!white}
            \Tile{4}{1}{red!25!white}
        \end{tikzpicture}
        && \begin{tikzpicture} 
            \def\s{0.4};
            \Tile{0}{3}{blue!25!white}
            \Tile{3}{1}{red!25!white}
            \Tile{4}{1}{red!25!white}
        \end{tikzpicture}\\
    \end{aligned}\end{displaymath}
\end{example}

\section{A bijection between $(m,n)$-words and two-toned tilings}
The motivation for this article is the observation that, for $n\geq 1$, the sets $\Words(m,n)$ and $\Tilings(m,n)$ have the same cardinality; see Lemmas~\ref{lem:mn_words_cardinality}~and~\ref{lem:two_toned_tilings_cardinality}.

\subsection{From $(m,n)$-words to two-toned tilings}
Let $\wb\in\Words(m,n)$. By definition, $\wb$ can be uniquely written as 
\begin{equation}\label{eq:mn_word_decomposition}
    \wb=w_1a^{(1)}w_2a^{(2)}w_3a^{(3)}\ldots a^{(k-1)}w_k a^{(k)}, 
\end{equation}    
where $k\in[n]$, each $a^{(i)}$ is a (possibly empty) sequence of $(m+1)$'s, the sum of the lengths of all $a^{(i)}$'s is $n-k$ and $w_1w_2\ldots w_k$ is a topless $(m,k)$-word. Let us write $\ell_i$ for the length of $a^{(i)}$.

The decomposition \eqref{eq:mn_word_decomposition} gives rise to a two-toned tiling $T_{\wb}$ as follows. 

\begin{construction}\label{constr:word_to_tiling}
    Let $\wb\in\Words(m,n)$ be decomposed as described in \eqref{eq:mn_word_decomposition}. Let $i\in[k]$ and set $w_{k+1}=0$.
    We define
    \begin{itemize}
        \item a blue strip $B_{i}$ of length $\ell_i+1$;
        \item a red strip $\hat{R}_{i}$ of length $w_i-w_{i+1}$.
    \end{itemize}
    Moreover, let $\hat{R}_{0}$ be a red strip of length $m-w_1$. The associated two-toned tiling $T_{\wb}$ is then derived from the sequence $\hat{R}_{0}B_{1}\hat{R}_{1}B_{2}\hat{R}_{2}\ldots B_{k}\hat{R}_{k}$ by replacing each red strip of length $s$ by a sequence of $s$ red squares.
\end{construction}

\begin{lemma}\label{lem:words_to_tilings_map}
    For $\wb\in\Words(m,n)$, the tiling $T_\wb$ is in $\Tilings(m,n)$.
\end{lemma}
\begin{proof}
    Since $w_1w_2\ldots w_k$ is a topless $(m,k)$-word it is guaranteed that $w_i-w_{i+1}\geq 0$ for all $i\in[k]$. Therefore, it follows that $T_\wb$ is a tiling using red squares and blue strips.
    
    To prove the claim, it thus remains to show that the number of red squares is $m$ and the length of the blue strips is $n$. But this follows immediately from the construction, because:
    \begin{itemize}
        \item the number of red squares is 
        \begin{displaymath}
            \lvert\hat{R}_0\rvert + \lvert\hat{R}_1\rvert + \cdots + \lvert\hat{R}_{k}\rvert 
                = m-w_1 + \sum_{i=1}^{k}\bigl(w_i-w_{i+1}\bigr) = m;
        \end{displaymath}
        \item the cumulated length of all blue strips is 
        \begin{displaymath}
            \sum_{i=1}^{k}\bigl(\ell_i+1\bigr) = k + \sum_{i=1}^{k}\ell_i = k+(n-k) = n.\qedhere
        \end{displaymath}
    \end{itemize}
\end{proof}

\begin{proposition}
    The map $\xi\colon\Words(m,n)\to\Tilings(m,n), \wb\mapsto T_\wb$ is a bijection.
\end{proposition}
\begin{proof}
    By Lemma~\ref{lem:words_to_tilings_map} and the uniqueness of the decomposition \eqref{eq:mn_word_decomposition}, the map $\xi$ is a well-defined map from $\Words(m,n)$ to $\Tilings(m,n)$. Moreover, by Construction~\ref{constr:word_to_tiling} this map is clearly injective. Now, since Lemmas~\ref{lem:mn_words_cardinality}~and~\ref{lem:two_toned_tilings_cardinality} state that the sets $\Words(m,n)$ and $\Tilings(m,n)$ have the same cardinality, this map must be a bijection.
\end{proof}

\begin{example}
    Consider the $(8,12)$-word $\wb=779329919900$. The decomposition \eqref{eq:mn_word_decomposition} is determined by the following values, where $\varepsilon$ denotes the empty word.
    \begin{center}\begin{tabular}{c|ccccccc}
        $i$ & $1$ & $2$ & $3$ & $4$ & $5$ & $6$ & $7$ \\
        \hline
        $w_i$ & $7$ & $7$ & $3$ & $2$ & $1$ & $0$ & $0$\\
        \hline
        $a^{(i)}$ & $\varepsilon$ & $9$ & $\varepsilon$ & $99$ & $99$ & $\varepsilon$ & $\varepsilon$\\
        \hline
        $\ell_i+1$ & $1$ & $2$ & $1$ & $3$ & $3$ & $1$ & $1$
    \end{tabular}\end{center}
    Then, the sequence $\hat{R}_{0}B_{1}\hat{R}_{1}B_{2}\hat{R}_{2}\ldots B_{7}\hat{R}_{7}$ induces the following two-toned tiling.
    \begin{center}\begin{tikzpicture}
        \def\s{.4};
        \Tile{1}{1}{red!25!white}
        \Tile{2}{1}{blue!25!white}
        \Tile{3}{2}{blue!25!white}
        \Tile{5}{1}{red!25!white}
        \Tile{6}{1}{red!25!white}
        \Tile{7}{1}{red!25!white}
        \Tile{8}{1}{red!25!white}
        \Tile{9}{1}{blue!25!white}
        \Tile{10}{1}{red!25!white}
        \Tile{11}{3}{blue!25!white}
        \Tile{14}{1}{red!25!white}
        \Tile{15}{3}{blue!25!white}
        \Tile{18}{1}{red!25!white}
        \Tile{19}{1}{blue!25!white}
        \Tile{20}{1}{blue!25!white}
    \end{tikzpicture}\end{center}
\end{example}

\subsection{From two-toned tilings to $(m,n)$-words}

Let us now explicitly describe the inverse map of $\xi$. 

\begin{construction}\label{constr:tiling_to_word}
    Let $T\in\Tilings(m,n)$, and let $B_1,B_2,\ldots,B_k$ denote its blue strips in order. Let $r_0$ denote the number of red squares before $B_1$ and for $i\in[k-1]$, let $r_i$ denote the number of red squares between $B_i$ and $B_{i+1}$. Since the total number of red squares is $m$, it follows that there must be $m-r_0-r_1-\cdots-r_{k-1}$ red squares after $B_k$.
    
    Let $\wb_T\defs w_1a^{(1)}w_2a^{(2)}\ldots w_ka^{(k)}$, where 
    \begin{align*}
        w_i & \defs m-\sum_{j=0}^{i-1}r_i,\\
        a^{(i)} & \defs \underbrace{(m+1)(m+1)\ldots(m+1)}_{|B_i|-1\;\text{times}}.
    \end{align*}
\end{construction}

\begin{lemma}\label{lem:tilings_to_words_map}
    For $T\in\Tilings(m,n)$, the word $\wb_T$ is in $\Words(m,n)$.
\end{lemma}
\begin{proof}
    It is sufficient to show that $w_1w_2\ldots w_k$ is a topless $(m,k)$-word, and that the total number of letters in $\wb_T$ is $n$. It follows immediately from the construction that $w_i\leq m$ for all $i$ and that $w_1\geq w_2\geq\cdots\geq w_k$ which establishes the fact that $w_1w_2\ldots w_k$ is a topless $(m,k)$-word. 

    For $i\in[k]$, let $\ell_i$ denote the number of copies of $m+1$ that are contained in $a^{(k)}$. Then, it follows that the number of letters of $\wb_T$ is 
    \begin{displaymath}
        k + \sum_{i=1}^{k}\ell_i = k + \sum_{i=1}^{k}\bigl(\lvert B_i\rvert-1\bigr) = \sum_{i=1}^{k}\lvert B_i\rvert,
    \end{displaymath}
    i.\;e., it equals the sum of the lengths of the blue strips. Since $T\in\Tilings(m,n)$, this number is exactly $n$.

    Therefore, $\wb_T\in\Words(m,n)$.
\end{proof}

\begin{proposition}
    The map $\xi^{-1}\colon\Tilings(m,n)\to\Words(m,n),T\mapsto\wb_T$ is a bijection.
\end{proposition}
\begin{proof}
    By Lemma~\ref{lem:tilings_to_words_map}, the map $\xi^{-1}$ is a well-defined map from $\Tilings(m,n)$ to $\Words(m,n)$. Moreover, Construction~\ref{constr:tiling_to_word} implies that this map is injective. Once again, Lemmas~\ref{lem:mn_words_cardinality} and \ref{lem:two_toned_tilings_cardinality} state that both sets $\Tilings(m,n)$ and $\Words(m,n)$ have the same cardinality, which proves the claim.
\end{proof}

\begin{example}
    Consider the following two-toned tiling $T$ of $6+11$:
    \begin{center}\begin{tikzpicture}
        \def\s{.4};
        \Tile{1}{1}{red!25!white}
        \Tile{2}{2}{blue!25!white}
        \Tile{4}{2}{blue!25!white}
        \Tile{6}{1}{red!25!white}
        \Tile{7}{1}{red!25!white}
        \Tile{8}{1}{blue!25!white}
        \Tile{9}{1}{red!25!white}
        \Tile{10}{3}{blue!25!white}
        \Tile{13}{1}{red!25!white}
        \Tile{14}{3}{blue!25!white}
        \Tile{17}{1}{red!25!white}
    \end{tikzpicture}\end{center}

    We get $r_0=1$, $r_1=0$, $r_2=2$, $r_3=1$, $r_4=1$, $r_5=1$. The lengths of the blue strips are $\lvert B_1\rvert=2$, $\lvert B_2\rvert=2$, $\lvert B_3\rvert=1$, $\lvert B_4\rvert=3$, $\lvert B_5\rvert=3$.

    Thus, we get $w_1=5$, $w_2=5$, $w_3=3$, $w_4=2$, $w_5=1$ so that 
    \begin{displaymath}
        \wb_T = 57573277177 \in\Words(6,11).
    \end{displaymath}
\end{example}

\section{Possible next steps}
In \cite{muehle24combinatorics,pilaud25hochschild}, the set of $(m,n)$-words was studied from an order-theoretic and geometric perspective. In particular it was shown that the set of $(m,n)$-words under componentwise order is a semidistributive lattice. This implies that the set of $(m,n)$-words admits a secondary order structure, the core label order as defined in~\cite{muehle19the}. 

A natural next step would be to transfer the order structure from the $(m,n)$-word lattice to two-toned tilings and investigate if the combinatorics of two-toned tilings helps with the understanding of the core label order of the $(m,n)$-word lattice.

\section*{Acknowledgements}
I thank Tom Copeland for observing computationally that the coefficient of $x^{n}$ of $\left(\frac{1-x}{1-2x}\right)^{m+1}$ enumerates the set of $(m,n)$-words and for providing multiple combinatorial explanations of this coefficients, in particular for pointing out the connection to two-toned tilings.


\end{document}